\documentclass[12pt]{article}
\usepackage{mathrsfs}

\usepackage{amsmath}
\usepackage{amssymb}
\usepackage{amsfonts}
\usepackage{graphicx}
\usepackage{subfigure}
\RequirePackage[colorlinks,citecolor=blue,urlcolor=blue]{hyperref}

\newtheorem{theorem}{Theorem}
\newtheorem{definit}{Definition}
\newtheorem{lemma}{Lemma}
\newtheorem{prop}{Proposition}

\newtheorem{remark}{Remark}

\begin{document}

\title{Random conformal welding for finitely connected regions}

\author{Shi-Yi Lan\footnote{Shi-Yi Lan was partially supported by the NSF of China (11161004) and NSF of Guangxi (2013GXNSFAA019015). E-mail address: shiyilan05@sina.com}
 and Wang Zhou\footnote{Wang Zhou was partially supported by a grant
R-155-000-151-112 at the National University of Singapore. E-mail address: stazw@nus.edu.sg}\\
{}\\
{\it Guangxi University for
Nationalities, and National University of Singapore}
}

\date{ }
\maketitle

\textbf{Abstract.} Given a finitely connected region $\Omega$ of the Riemann sphere  whose complement
consists of $m$ mutually disjoint closed disks $\bar{U}_j$,  the random homeomorphism $h_j$ on the boundary component $\partial U_j$
 is constructed using the exponential Gaussian free field. The existence and uniqueness of random conformal welding of $\Omega$
with $h_j$ is established by investigating a non-uniformly elliptic Betrami equation with a random complex
dilatation. This generalizes the result of Astala, Jones, Kupiainen
and Saksman to  multiply connected domains.

\textbf{Keywords:} Random welding;
quasiconformal mapping; Gaussian free field; SLE.

\textbf{MSC(2010):} 30C62, 60D05.


\section{Introduction}

Over the last decades there has been great interest in conformally invariant fractals which could arise as scaling
limits of discrete random processes in the complex plane or the Riemann sphere.  One of the most important
examples of confromally invariant fractals is the Schramm Loewner Evolution (SLE) introduced by
Schramm \cite{sc} in 2000. There are several different
versions of SLE, among which chordal SLE and radial SLE are the most well-known.
A chordal SLE trace describes a
random curve evolving in the upper plane from one point on the boundary to another point on the boundary.
A radial SLE trace describes a random curve evolving in the disk from a point on the boundary to
an interior point. The behavior of the SLE trace depends on a real parameter $\kappa>0$. If $\kappa\in (0,4]$, the trace is a simple curve; if $\kappa\in (4,8)$, the trace has self-intersections; and if $\kappa\in [8,\infty)$,
the trace is space-filling. For more information on $\mbox{SLE}_\kappa$
and related topics, see \cite{law,ros} and \cite{kan} etc..

Moreover, many two-dimensional random lattice paths from statistical physics have been proved
to have $\mbox{SLE}_\kappa$ curves as their scaling limits when the mesh of the grid tends
 to $0$, such as the critical site percolation exploration path \cite{sm1,sm1b, CN07}, loop erased
random walks and uniform spanning tree Peano paths \cite{lsw}, the harmonic explorer’s
path \cite{scs1}, the  chordal contour lines of the discrete Gaussian free field \cite{scs2}, the
interfaces of the FK Ising model \cite{sm2}.

Besides conformal invariance, the SLE paths evolve in the domain. In other words, the SLE paths are indexed by capacity or natural parametrization \cite{ls, lz}. However there are other conformally invariant curves which are independent of an auxiliary time.
Recenly, Astala, Jones, Kupiainen and Saksman \cite{ajks} used  the idea of conformal welding
to construct a random family of closed conformally invariant curves in the Riemann sphere based
on a method of Lehto \cite{le,aim} and a result \cite{js} on the conformal
removability of H\"{o}lder curves. This family of random curves obtained in \cite{ajks} is stationary
for each inverse temperature less than a certain critical value.
Instead of the white noise representation for the Gaussian free field, in a similar manner
Tecu \cite{te} extended the work of
Astala, Jones, Kupiainen and Saksman \cite{ajks} to the situation of criticality using a vaguelet representation
of the Gaussian free field. In the meantime, we also note that Sheffield \cite{she1} investigated a conformal welding of two Liouville quantum gravity random surfaces
using a totally different approach. Two quantum surfaces of normalized quantum area are welded together by
matching quantum length on the boundaries. This results in an SLE interface.

However, these results mentioned above deal  with random conformal weldings only for simply connected domains.
In this paper we are concerned with the
random conformal welding associated with finitely connected domains.
Fix a positive integer $m\geq 2$, suppose that
$\Omega$ is an $m$-connected region of the Riemann sphere $\mathbb{S}^2=\mathbb{C}\cup\{\infty\}$
whose complement $\mathbb{S}^2\backslash \Omega$ is a union of $m$ disjoint closed disks $\bar{U}_j, j=1,2,\dots,m$.
Our main goal is to establish the existence and uniqueness theorem for the random conformal
welding of $\Omega$, which will yield $m$ mutually disjoint  random Jordan curves on $\mathbb{S}^2$.
Our method involves modifications and generalizations of those in \cite{ajks}.
We first construct random measures by limiting processes via the exponential Gauss free fields
restricted to the boundary components $\partial U_j$
of $\Omega$, which allow us to define random homeomorphisms $\psi_j$ on $\partial U_j$.
Next we will solve the conformal welding problem
of $\Omega$ with  $\psi_j$, i.e., to seek a random conformal mapping $f$ from $\Omega$ into
$\mathbb{S}^2$ and random conformal mapping $g_j$ from $U_j$ into $\mathbb{S}^2$ such that
$g_j(\zeta)=f\circ \psi_j(\zeta)$  when $\zeta\in \partial U_j, j=1,2,\dots,m$, and the images
$\Gamma_j=f\circ \psi_j(\partial U_j)=g_j(\partial U_j)$ are the desired random Jordan curves.

To this end, applying the technique
of Beurling-Ahlfors extension \cite{bea} we extend the random homeomorphisms $\psi_j:\partial
U_j\rightarrow\partial U_j$ to $\Omega$,
which leads to a quasiconformal homeomorphism of $\mathbb{S}^2$
with a random complex dilatation $\lambda$. Thus
the welding problem is reduced to solving a non-uniformly elliptic Beltrami equation with
$\lambda$. We prove the existence and uniqueness of the solution to
the Beltrami equation using the techniques of Lehto \cite{le,aim}, and
the conformal removability for boundary components of multiply connected domain.
Hence we get the existence and uniqueness of random conformal welding
associated with $\Omega$, which generalizes the result of \cite{ajks} to finitely connected domains. Our main result can be summarized as follows.

{\it With probability one, there exist random conformal mappings
$f:\Omega\rightarrow\mathbb{S}^2$ and
$g_j:U_j\rightarrow\mathbb{S}^2$, respectively, such that their boundary values satisfy
$g_j(\zeta)=f\circ \psi_j(\zeta)$  when $\zeta\in \partial U_j, j=1,2,\dots,m$, which produce $m$ mutually disjoint random Jordan curves $\Gamma_j=f\circ \psi_j(\partial U_j)=g_j(\partial U_j)$. Moreover, almost surely in $\omega$, these Jordan curves $\Gamma_j$ are unique up to a M\"{o}bius transformation
$\chi=\chi_\omega$ of the Riemann sphere $\mathbb{S}^2$.}

We refer to Theorem \ref{thm2} in Section \ref{RCW} for a complete statement involving Gaussian free fields.

Although our result is a generalization of  \cite{ajks,te}, there is a big difference between our paper and \cite{ajks,te}, which can be summarized by the following three points.
\begin{enumerate}
\item
The random conformal welding
for simply connected regions is discussed in \cite{ajks,te}, which generates a random
Jordan curve in $ \mathbb{C}$.
Instead, we here deal with the random conformal welding associated with the finitely connected domain $\Omega$, which leads to  $m$ mutually disjoint  random  Jordan curves in $ \mathbb{S}^2$.
\item
The construction of the extension mapping in \cite{ajks,te}
involves only one random homeomorphism,
while the corresponding extension mapping in this paper is produced
by  $m$ random homeomorphisms.
We also apply the fact that these random homeomorphisms are invariant
under conformal transformations of $\mathbb{S}^2$.
\item
The unique solution to the Beltrami equation in \cite{ajks,te}
is determined by the conformal removability of a H\"{o}lder Jordan curve,
whereas the uniqueness of the corresponding solution in the present paper is
obtained by showing that the
$m$ boundary components of the multiply connected domain are conformally removable. To achieve this we appeal to a result of \cite{mar,doy}, i.e., a conformal mapping from
a finitely connected domain can be expressed as a
composition of finitely many conformal mappings of simply connected domains.
\end{enumerate}

In addition, it is
worth to point out that the deterministic version of conformal welding for the
finitely connected region $\Omega$ has been discussed by Marshall \cite{mar}  using the
geodesic zipper algorithm and Koebe's iterative method to compute conformal mappings.
In a sense the current paper may be viewed as a random version of \cite{mar}.

This paper is organized as follows. In Section \ref{cwgff} we shall briefly introduce the definitions
of conformal welding and  Gaussian free field, and state some useful results. In Section \ref{crm}
we construct random homeomorphisms on the boundary components $\partial \bar{U}_j$ of
the multiply connected domain $\Omega$
by means of  the exponential Gaussian free field. In Section \ref{erh} we show that there exists a solution to the
non-uniformly elliptic Beltrami equation with a random complex dilatation $\lambda$. Moreover, this solution
is unique up to  M\"{o}bius transformations of $\mathbb{S}^2$.
In Section \ref{RCW} we establish random conformal welding theorems
for the multiply connected domain $\Omega$, which generates $m$ random conformal welding
curves $\Gamma_j$ on the Riemann
sphere $\mathbb{S}^2$.

\section{ Conformal weldings and Gaussian free field} \label{cwgff}
In this section, we will briefly review  some basic concepts
related to conformal welding  and Gaussian free field, and provide some useful results;
see \cite{lv,ha2,ras,bis,aim,she,scs2} for more details.

\subsection{Conformal weldings}

Consider the Riemann sphere $\mathbb{S}^2=\mathbb{C}\cup\{\infty\}$. The conformal welding arises usually in the following two cases.

The first case is the \textbf{conformal welding of simple connected domains.}
Given two disjoint Jordan domains $\Omega_1$ and $\Omega_2$ in $\mathbb{S}^2$
and a homeomorphism $\psi:\partial\Omega_2\rightarrow\partial\Omega_1$,
one can attach $\Omega_1$ and $\Omega_2$ by identifying points $\zeta$ of $\partial\Omega_2$ with points $\psi(\zeta)$ of $\partial\Omega_1$.
The mapping $\psi$ is called a welding homeomorphism if there exist conformal mappings $f$ and $g$ of $\Omega_1$ and $\Omega_2$, respectively,
onto complementary regions of $\mathbb{S}^2$ such that
\begin{equation}\label{eq1} g(\zeta)=f\circ\psi(\zeta)\end{equation}
for each $\zeta\in \partial \Omega_2$.
If $\psi$ is a welding homeomorphism, then it induces a conformal welding.
Since the two regions push and pull against
one another as they find their new positions, this yields a Jordan curve $\Gamma$ known as the
conformal welding curve.

Conversely, if we are given a homeomorphism $\psi$ which maps $\partial\Omega_2$ onto
$\partial\Omega_1$, the conformal welding problem is to seek a conformal welding curve $\Gamma$ and
conformal mappings $f:\Omega_1\rightarrow\Omega_1^*$ and $g:\Omega_2\rightarrow\Omega_2^*$
such that their boundary values satisfy (\ref{eq1}), where $\Omega_1^*\cup \Gamma\cup\Omega_2^*=\mathbb{S}^2$.
The Conformal Welding Theorem \cite{lv,gl,ha2} tells us that
if $\Omega_1$ and $\Omega_2$ are both disks and $\psi$ is
a quasi-symmetric mapping, then a conformal welding will exist and the conformal welding curve will be a quasicircle.

\begin{remark}
  Williams \cite{wi2} constructed discrete conformal weldings of the first case
which converge uniformly on compact subsets to their continuous counterparts.
The corresponding random conformal weldings were investigated by
Astala, Jones, Kupiainen and Saksman \cite{ajks}, which produces a conformally invariant random family of
closed curves in $\mathbb{S}^2$. Tecu \cite{te} extended the result of \cite{ajks}
to the situation of criticality.
\end{remark}

The second case is the \textbf{conformal welding of finitely connected regions.}
Fix a natural number $m\geq 2$, and suppose we are given $m$ disjoint
simply connected domains $\Omega_1,\Omega_2,\dots,\Omega_m$ on $\mathbb{S}^2$,
which are bounded by disjoint Jordan curves.
The Riemann mapping Theorem and Koebe's Theorem imply that there exist conformal mappings $g_j$ of disks $U_j$
onto $\Omega_j$ and a conformal mapping $f$ of $\mathbb{S}^2\setminus \cup_{j=1}^m \bar{U}_j$
onto $\mathbb{S}^2\setminus \cup_{j=1}^m\bar{\Omega}_j$.
We call the mappings $\psi_j=f^{-1}\circ g_j:\partial U_j\rightarrow\partial U_j$ the welding homeomorphisms
associated with the multiply connected domain $\mathbb{S}^2\setminus \cup_{j=1}^m \bar{U}_j$ for $j=1,2,\dots,m$.

On the other hand, consider an $m$-connected region
$\Omega=\mathbb{S}^2\backslash \cup_{j=1}^m \bar{U}_j$.
If we are given $m$ homeomorphisms $\psi_j:\partial U_j\rightarrow\partial U_j$
for $j=1,2,\dots,m$, then the conformal welding problem for $\Omega$ and $U_j$
is to find a conformal mapping $f$ from $\Omega$ into
$\mathbb{S}^2$ and conformal mapping $g_j$ from $U_j$ into $\mathbb{S}^2$ such that
\begin{equation}\label{eq2}
g_j(\zeta)=f\circ\psi_j(\zeta)
\end{equation}
for all $\zeta\in \partial U_j, j=1,2,\dots,m$.  Consequently, this will lead to
a sphere with patches $\Omega_j$, bounded by Jordan curves $L_j=f(\partial U_j)=g_j(\partial U_j)$
for $j=1,2,\dots,m$; see Figure 1.

\begin{figure}
\centering
\includegraphics[width=5.5in]{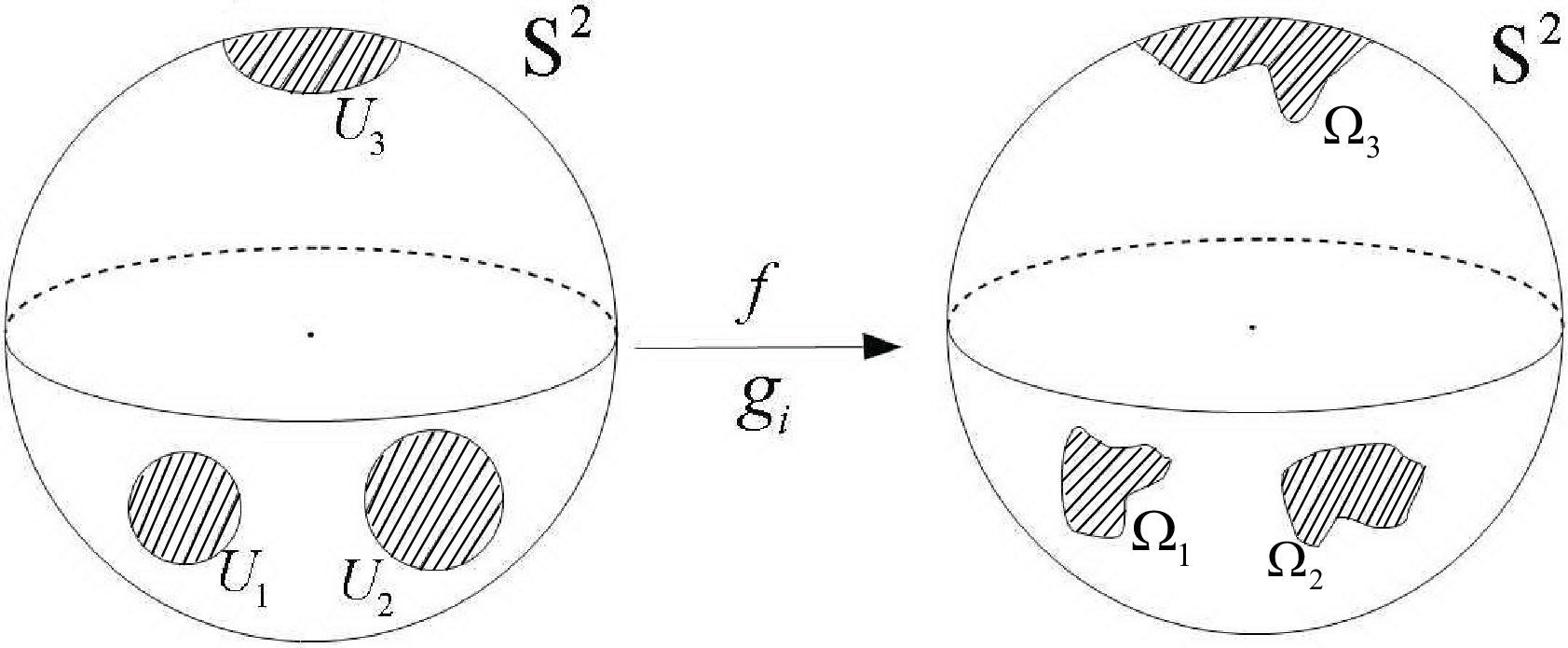}\\
\centerline{\small\rm} Figure 1. A conformal welding exists for $\psi_j$ if there are conformal
mappings $f$ and $g_j$ onto complementary regions of $\mathbb{S}^2$ whose boundary values satisfy
$g_j=f\circ \psi_j$ for $j=1,2,3$.
\end{figure}

\

\begin{remark} The numerical implementation of the second case, where $U_j$ are
all equal to the unit disk for $j=1,2\dots,m$, was discussed by Marshall \cite{mar}  using the
geodesic zipper algorithm and Koebe's iterative method to compute conformal mappings.
In the present paper our aim is to establish a random version of conformal welding theorem
for the second case, which results in $m$ mutually disjoint random Jordan curves in $\mathbb{S}^2$.
\end{remark}

\subsection{Gaussian free field}

 The two-dimensional Gaussian free field (GFF) is a two-dimensional-time analog of Brownian motion,
which may be viewed as a Gaussian random variable on an infinite dimensional space.

\begin{definit}
\label{Definition1} For a given planar domain $D\subset\mathbb{S}^2=\mathbb{C}\cup\{\infty\}$ let $H_s(D)$ be  the space
of $C^\infty$ real-valued functions with compact support in $D$, and let $H(D)$ be its Hilbert
space completion under the Dirichlet inner product
\[(f_1,f_2)_{\nabla}:=(2\pi)^{-1}\int_D\nabla f_1(z)\cdot\nabla f_2(z)dz,\]
where $dz$ refers to area measure. We define an instance $h$ of the Gaussian free field (GFF) to be the
formal sum
\begin{equation}\label{eq1a}h=\sum_{j=1}^\infty\alpha_jf_j,\end{equation}
where the $\alpha_j$ are i.i.d. one-dimensional standard (unit variance, zero mean) Gaussian random variables, and $\{f_j, j=1,2,\cdots\}$ is any orthonormal
base of $H(D)$.
\end{definit}
The sum \eqref{eq1a} does not converge within $H(D)$ almost surely, since $\sum|\alpha_j|^2$ is
infinite almost surely. However, it does converge almost surely in the space of distributions on $D$. That is,
the limit $(\sum_{j=1}^\infty\alpha_jf_j, g)$ almost surely exists for all $g\in H_s(D)$, and the limit
value as a function of $g$ is almost surely continuous on $H_s(D)$. For each $f\in H(D)$, $(h,f)_\nabla$
is a mean zero  Gaussian random variable, and
\begin{equation}\label{inner}
\mbox{Cov}((h,f_1)_\nabla,(h,f_2)_\nabla)=(f_1,f_2)_\nabla
\end{equation}
for any $f_1, f_2\in H(D)$. The collection of random variables $(h,f)_\nabla$ for $f\in H(D)$ is thus a Hilbert space with the inner product \eqref{inner}.

\textbf{Conformal invariance.} Let $\phi$ be a conformal mapping from $D$ to another planar domain $\tilde{D}$. Then an elementary
change of variables calculation gives that
\[\int_{\tilde{D}}\nabla(f_1\circ\phi^{-1})\cdot\nabla(f_2\circ\phi^{-1})dw=\int_D(\nabla f_1\cdot\nabla f_2)dz\]
for any $f_1,f_2\in H_s(D)$. Taking the completion to $H(D)$, we see that the Dirichlet inner product is
invariant under conformal transformations of $D$. This implies that the two-dimensional GFF possesses the
conformal invariance property.

\textbf{Representation of covariance.} For a fixed $\zeta\in D$, let $\tilde{G}_\zeta(z)$ be the harmonic extension to $z\in D$ of the function of
$z$ on $\partial D$ given by $-\log|z-\zeta|$. Then Green's function in the domain $D$ is defined by
\[G(\zeta,z)=-\log|z-\zeta|-\tilde{G}_\zeta(z).\]
It is known that if $\zeta\in D$ is fixed, then Green's function $G(\zeta,z)$ may be viewed
as a distributional solution of the Poisson equation $\Delta G(\zeta,\cdot)=-2\pi\delta_\zeta(\cdot)$
with zero boundary conditions.

For each $g\in H_s(D)$,  we define a function $\Delta^{-1}g$ on $D$ by
\[\Delta^{-1}g(\cdot):=-\frac{1}{2\pi}\int_D G(\cdot,z)g(z)dz.\]
Then $\Delta^{-1}g(\cdot)$ is a $C^\infty$ function in $D$ whose Laplacian is $g$.
If $f_1=-\Delta^{-1}g_1$ and $f_2=-\Delta^{-1}g_2$, then integration by parts gives that
\[(f_1,f_2)_\nabla=\frac{1}{2\pi}(g_1,-\Delta^{-1}g_2),\]
where $(\cdot,\cdot)$ denotes the standard inner product in $L^2(D)$. Note that each $h\in H(D)$
is naturally a distribution, since one may define the map $(h,\cdot)$ by $(h,g):=2\pi(h,-\Delta^{-1}g)_\nabla$
for any $g\in H_s(D)$.
It is easy to see that $-\Delta^{-1}g\in H(D)$. Thus, if $-\Delta f=g$, then we may write $(h,g)=2\pi(h,f)_\nabla$.
This implies
\[\mbox{Cov}((h,g_1),(h,g_2))=(2\pi)^2(f_1,f_2)_\nabla.\]
Hence, it follows that
\[\mbox{Cov}((h,g_1),(h,g_2))=\int_{D\times D}g_1(\zeta)G(\zeta,z)g_2(z)d\zeta dz.\]
So $G(\zeta,z)$ is also  the integral kernel of covariance $\mbox{Cov}((h,f_1),(h,f_2))_\nabla$ for any
$f_1,f_2\in H(D)$.

Notice that in this paper we consider only the restriction of the Gaussian free field (\ref{eq1a}) to
the boundary components $\partial U_j$, $j=1,\cdots,m$, of the multiply connected region $\Omega$.

\section{The construction of random homeomorphisms} \label{crm}

In this section we will describe how to construct random
homeomorphisms on the boundary components of the finitely connected domain $\Omega$ in $\mathbb{S}^2$.
First, the restriction of Gaussian free field (\ref{eq1a}) $(D=\mathbb{C})$ to each boundary component of $\Omega$ may be given by
a concrete expression, and the latter can be further expressed in terms of a white noise representation.
Next we use the white noise representations of (\ref{eq1a})
to construct random measures on the boundary components of $\Omega$,
which can be viewed as martingale limits
of products of exponentials of independent Gaussian fields. Thus we may define random
homeomorphisms on the boundary components of $\Omega$ and derive some useful results. In
particular, we show that these random homeomorphisms are conformally invariant. This generalizes
the random homeomorphism of the unit circle constructed in \cite{ajks} to the case of $m$
mutually disjoint circles with finite radii in $\mathbb{C}$.

\subsection{The representations of GFF on boundary components}

As before, let $\Omega\subset\mathbb{S}^2$ be an $m$-connected domain
whose complement is a union of disjoint closed disks $\bar{U}_j, j=1,2,\dots,m,$ for any fixed integer $m\geq 2$,
and assume that $\infty$ belongs to $U_m$.
For convenience, we will work on the complex plane $\mathbb{C}$
instead of the Riemann sphere $\mathbb{S}^2$, keeping in mind that $\infty$ corresponds to one point
on $\mathbb{S}^2$. This implies that $U_j, j=1,2,\dots,m-1,$ are bounded,
and that $U_m$ are unbounded in $\mathbb{C}$ where $\infty$ may be viewed as a point
in $U_m$. Thus we may write $U_j=\{z\in\mathbb{C}:|z-a_j|<r_j\}$
where $|a_j|<\infty,0<r_j<\infty$ for $j=1,2,\dots,m-1$, and assume that
the exterior of $U_m$ is equal to the disk $U_m^c=\{z\in\mathbb{C}:|z-a_m|<r_m\}$
where $|a_m|<\infty,0<r_m<\infty$.
Hence  $\Omega$ can be written as $\Omega=U_m^c\backslash \cup_{j=1}^{m-1} \bar{U}_j$.

We first give the concrete representations of traces of GFF on the boundary components $\partial U_j$
of $\Omega$. Set $h_j=h|_{\partial U_j}, j=1,2,\dots,m$, which
may be viewed as the restriction of the 2-dimensional GFF on $\mathbb{C}$ to $\partial U_j$.
The covariance functions of $h_j$, $j=1,2,\dots,m$,
have the integral kernels
\begin{equation}\label{equ1b}G_{h_j}(\zeta,\zeta')=-\log|\zeta-\zeta'|,\quad \zeta,\zeta'\in\partial U_j\end{equation}
If we let $\zeta=a_j+r_je^{i2\pi t},\zeta'=a_j+r_je^{i2\pi t'}, t,t'\in [0,1)$,
where $a_j$ and $r_j$ are the centers and the radii of $U_j, j=1,2,\dots,m-1,$ respectively, and $a_m$ and $r_m$ are
the center and the radius of $U_m^c$
respectively,
then it follows from (\ref{equ1b}) that the covariance functions of $h_j$, $j=1,2,\dots,m$, may take the forms
\begin{equation}\label{equ1}
G_{h_j}(t,t')=-\log2r_j\sin\pi|t-t'|,~t,t'\in[0,1)
\end{equation}
when $\partial U_j$ is identified with $[0,1)$.
A direct computation gives that the fields $h_j$, $j=1,2,\dots,m$, with covariances (\ref{equ1})
can be expressed by the Fourier expansions
\begin{equation}\label{equ2}
h_j=\sum_{n=1}^\infty\frac{r_j^n}{\sqrt{n}}(\alpha_n^{(j)}\cos2\pi nt+\beta_n^{(j)}\sin2\pi nt),
t\in[0,1)\end{equation}
where $\alpha_n^{(j)},\beta_n^{(j)}\sim N(0,1),n\geq 1$ are independent standard
Gaussian random variables. We remark that if $\partial U_j$ are all equal to the unit circle for $j=1,2,\dots,m$,
i.e., $a_j=0,r_j=1$, then (\ref{equ1})
and (\ref{equ2}) will become the resutls discussed in \cite{ajks}.

Next, we will describe further that the formula (\ref{equ2}) can be expressed
by white noise representations.  A  white noise $Y$ in the upper half-plane $\mathbb{H}$ is a centered
Gaussian process, indexed by sets with finite hyperbolic area measure in $\mathbb{H}$,
whose covariance structure is given by the hyperbolic area measure of the intersection of sets.
We will need a periodic version of $Y$, which can be identified with a white noise
on $[0,1)\times\mathbb{R}_+$. To be more specific, set
\[W=\{(x,y)\in\mathbb{H}:-\frac{1}{2}<x<\frac{1}{2}~
\mbox{and}~y>\frac{2}{\pi}\tan|\pi x|\}\]
and
\[V=\{(x,y)\in\mathbb{H}:-\frac{1}{4}<x<\frac{1}{4}~\mbox{and}~ 2|x|<y<\frac{1}{2}\}.\]
For a small positive number $\epsilon>0$,  we define two random fields $Y^\epsilon(x)$ and  $Z^\epsilon (x)$ by
\begin{equation}\label{equ2a}
Y^\epsilon(x):=Y(x+W_\epsilon),~x\in [0,1)
\end{equation}
where $W_\epsilon=\{(x,y)\in W:y>\epsilon\}$, and
\begin{equation}\label{equ2b} Z^\epsilon (x):=Y(x+V_\epsilon),~x\in [0,1)\end{equation}
where $V_\epsilon=\{(x,y)\in V:y>\epsilon\}$, respectively.
Then the covariance functions of two fields $Y_\epsilon(x)$ and  $Z_\epsilon (x)$ can be expressed by
\begin{equation}\label{equ2a1}\mathbb{E}(Y^\epsilon(x_1)Y^\epsilon(x_2))=\mu((x_1+W_\epsilon)\cap\bigcup
\limits_{z\in\mathbb{Z}}(x_2+W_\epsilon+n))\end{equation}
and
\begin{equation}\label{equ2a2}
\mathbb{E}(Z^\epsilon(x_1)Z^\epsilon(x_2))=\mu((x_1+V_\epsilon)\cap\bigcup\limits_{z\in\mathbb{Z}}(x_2+V_\epsilon+n))
\end{equation}
respectively, where $\mu$
denotes the hyperbolic area measure in $\mathbb{H}$, given by
$\mu(dxdy)=dxdy/y^2$.
Then we have the following lemma, which formulates that the restriction
of $h$ to $\partial U_j$ can be represented by the white noise.

\begin{lemma}
\label{lem1}
(i). For each $h_j=h|_{\partial U_j}, j=1,2,\dots,m$,
there exists a version $Y_j^\epsilon (x)$ of the white noise (\ref{equ2a})
which converges weakly to a random field
$Y_j$ as $\epsilon\rightarrow 0$ satisfying $Y_j\sim h_j+G_j$,
where $G_j\sim N(0,2r_j\log2)$ is a scalar Gaussian factor,
independent of $h_j$, and $r_j$ denotes the radius of circle $\partial U_j$.

(ii). There exists a version $Z_j^\epsilon(x)$ of the white noise (\ref{equ2b}) corresponding to
  $Y_j^\epsilon (x)$
such that
\[w_j:=\sup\limits_{x\in [0,1),\epsilon\in(0,1/2]}|Z_j^\epsilon(x)-Y_j^\epsilon(x)|<\infty \quad \mbox{a.~s.}.\]
Moreover, $\mathbb{E}e^{q w_j}<\infty$ for any $q>0$.
\end{lemma}
\noindent\textbf{Proof.} Since the Gaussian free field $h$ on $\mathbb{C}$
is conformally invariant,
the distribution of $h_j=h|_{\partial U_j}$ is identified with
one of $h|_{\partial \mathbb{U}}$ for each $\partial U_j$,
where $\partial\mathbb{U}$ denotes the unit circle. Thus
similar to the proof of \cite[Lemma 3.4]{ajks},
according to Duldey's theorem we may conclude (i). A straightforward
computation through (\ref{equ2a1}) replacing $Y^\epsilon$ and $W^\epsilon$ by $Y_j$
and $W$, respectively,
combined with (\ref{equ1}), gives
that $Y_j$ has the same distribution as
$h_j+G_j$.
Note that each $Y_j^\epsilon$ has the same law as $Y^\epsilon(x)$, while $Y^\epsilon(x)$
is equivalent to $H_\epsilon(x)$ in \cite{ajks}. Therefore in the light of (\ref{equ2a1}) and (\ref{equ2a2}),
following the proof of \cite[Lemma 3.5]{ajks}
we can conclude that (ii) holds.
 $\hfill\square$

\subsection{The homeomorphisms from random measures} \label{home}

The stationarity of $Y_j^\epsilon(x)$ in Lemma 1(i) implies that $\mbox{Var}(Y_j^\epsilon(x))$
is independent of $x\in [0,1)$, where $\mbox{Var}(Y_j^\epsilon(x))=\mathbb{E}|Y_j^\epsilon(x)|^2$.
Assume that $\beta_j>0, j=1,2,\dots,m$, and let $M([0,1))$ denote
the set of bounded Borel measures on $[0,1)$.
For any $x\in [0,1)$ and each $g\in M([0,1))$, consider
the following processes
\[X_j^\epsilon:=e^{\beta_j Y_j^\epsilon(x)-\beta_j^2\mbox{Var}(Y_j^\epsilon(x))/2}\]
and
\[\hat{X}_j^\epsilon:=\int_0^{1}e^{\beta_j Y_j^\epsilon(x)-\beta_j^2\mbox{Var}(Y_j^\epsilon(x))/2}
g(x)dx.\]
It is easy to see that
$X_j^\epsilon$ and $\hat{X}_j^\epsilon$ are $L^1$-martingales with respect to decreasing $\epsilon\in
(0,1/2]$. The martingale convergence theorem gives that $X_j^\epsilon$ and $\hat{X}_j^\epsilon$
converge almost surely as $\epsilon\rightarrow 0$, and that their $L^1$-norms stay bounded.
This, combined with Lemma 1(i),
gives rise to random measures $\tau_j$
on $[0,1)$, which can be defined  by
\begin{equation}\label{equ3} \tau_j(dx):=
\lim_{\epsilon\rightarrow 0}e^{\beta_j Y_j^\epsilon(x)-\beta_j^2\mbox{Var}(Y_j^\epsilon(x))/2}e^{-\beta_j G_j}\frac{dx}{2^{r_j\beta_j^2}}
\quad \mbox{w}^* \quad\mbox{in}\quad M([0,1))\quad\mbox{ a.~s.}, \end{equation}
where $G_j\sim N(0,2r_j\log 2)$, $j=1,2,\dots,m$, are Gaussian random variables. The limit measures
$\tau_j$, $j=1,2,\dots,m$, are weak*-measurable in the sense that  the
integrals $\int_0^1 g(x)\tau_j(dx)$ are well-defined random variables
for all $g\in C([0,1))$.

In addition, with the same reason as above we may define the random measure $\nu_j$
on $[0,1)$ corresponding to $Z_j^\epsilon(x)$
in Lemma 1(ii) by
\begin{equation}\label{equ3a}
\nu_j(dx):=\lim_{\epsilon\rightarrow 0}e^{\beta_jZ_j^\epsilon(x)-\beta_j^2\mbox{Var}(Z_j^\epsilon(x))/2}dx
\quad \mbox{w}^* \quad\mbox{in}\quad M([0,1))\quad\mbox{ a.~s.},
\end{equation}
where $\mbox{Var}(Z_j^\epsilon(x))=\mathbb{E}|Z_j^\epsilon(x)|^2$. Here is a lemma about the two measure $\tau_j$ and $\nu_j$.

\begin{lemma}
\label{lem2} (a) There exist almost surely  positive finite random variables $\mathcal{G}_j$, $j=1,2,\dots,m$,
satisfying $\mathbb{E}\mathcal{G}_j^q<\infty$  for any $q\in\mathbb{R}$, such that
\[\frac{1}{\mathcal{G}_j}\tau_j(B)\leq \nu_j(B)\leq \mathcal{G}_j\tau_j(B),\]
for any Borel set $B\subset [0,1)$.

(b) For each fixed $\beta_j<\sqrt{2}$, there exist $a_j^{(l)}=a_j^{(l)}(\beta_j), \tilde{a}_j^{(l)}=\tilde{a}_j^{(l)}
(\beta_j)>0, l=1,2,$
and almost surely finite random constants $c_j=c_j(\omega,\beta_j),\tilde{c}_j=\tilde{c}_j(\omega,\beta_j)>0$
such that
\[\frac{1}{c_j(\omega,\beta_j)}|I|^{a_j^{(1)}}\leq\tau_j(I)\leq c_j(\omega,\beta_j)|I|^{a_j^{(2)}},~
\frac{1}{\tilde{c}_j(\omega,\beta_j)}|I|^{\tilde{a}_j^{(1)}}\leq\nu_j(I)\leq \tilde{c}_j(\omega,\beta_j)|I|^{\tilde{a}_j^{(2)}}\]
for each subinterval $I\subset[0,1)$.
\end{lemma}

\noindent\textbf{Proof.} By Lemma \ref{lem1}(ii) and the stationary properties of the fields $Y_\epsilon(x)$ and $Z_\epsilon(x)$,
we deduce that (a) holds. It is easy to see from (\ref{equ3}) and (\ref{equ3a}), combined with the constructions
$Y_j^\epsilon$ and $Z_j^\epsilon$, that $\tau_j$ and $\nu_j$ are the same measures as $\tau$ and $\nu$ in \cite{ajks},
respectively.  So it follows from \cite[Theorem 3.7]{ajks} that (b) holds. $\hfill\square$

Now we are able to define the random homeomorphism on the boundary component
$\partial U_j$ of  $\Omega$, which is guaranteed by Lemma \ref{lem2}(b).

\begin{definit}
\label{Definition2} Let $a_j$ and $r_j$ be the center and radius of the circle $\partial U_j$,  and let $\beta_j<\sqrt{2}, j=1,2,\dots,m$.
Then the random homeomorphism
$\psi_j:\partial U_j\rightarrow\partial U_j$ can be obtained by setting
\begin{equation}\label{equ7}
\psi_j(a_j+r_je^{2\pi ix})=a_j+r_je^{2\pi ip_j(x)},
\end{equation}
where $p_j(x)$ is given by
\begin{equation}\label{equ8}
p_j(x)=p_{\beta_j}(x)=\frac{\tau_j([0,x])}{\tau_j([0,1))}
\end{equation}
for $x\in[0,1)$ and is extended periodically over the real line $\mathbb{R}$ for $j=1,2,\dots,m$.
\end{definit}

Thus we obtain
$m$ random homeomorphisms $\psi_j$, $j=1,2,\dots,m$, on the boundary components $\partial U_j$ of  $\Omega$,
which have the following properties.

\begin{lemma}
\label{lem3} Suppose that $\beta_j<\sqrt{2}, j=1,2,\dots,m$.
Then (i) almost surely both $\psi_j$ and its inverse mapping $\psi_j^{-1}$ are H\"{o}lder
continuous for any $1\leq j\leq m$;
(ii) the distribution of  $\psi_j$ is invariant under any M\"{o}bius transformation of
$\mathbb{S}^2$.
\end{lemma}

\noindent\textbf{Proof.}
It is obvious that (i) can follow from the definition of $\psi_j$ and Lemma \ref{lem2}(b).
In the following we will prove (ii), i.e., to show that $\chi\circ \psi_j\circ\chi^{-1}$ and $\psi_j$
have identical distributions for any M\"{o}bius transformation $\chi$ of
$\mathbb{S}^2$.
First, we show that $h_j$
is conformally invariant for $j=1,2,\dots,m$.  Indeed,
$h_j =h|_{\partial U_j}$  may be viewed as
the restriction of the 2-dimensional GFF $h_j$ on $\mathbb{C}$ to $\partial U_j$, whose
covariance function $G_{h_j}(\zeta,\zeta')$ is given by (\ref{equ1b}).

Let $\chi:\mathbb{S}^2\rightarrow \mathbb{S}^2$ be any  M\"{o}bius transformation,
and set $\tilde{h}_j=h|_{\partial \tilde{U}_j}$,
where $\partial \tilde{U}_j=\chi(\partial U_j)$, which could be also identified with $[0,1)$.
Then from the expressions (\ref{equ1b}) of $G_{h_j}(\zeta,\zeta')$ and
$G_{\tilde{h}_j}(\chi(\zeta),\chi(\zeta'))$
we deduce that
\[G_{\tilde{h}_j}(\chi(\zeta),\chi(\zeta'))=G_{h_j}(\zeta,\zeta')+Q_{j,1}(\zeta)+Q_{j,2}(\zeta'),\]
where $Q_{j,1}$ (respectively, $Q_{j,2}$) are independent of $\zeta'$ (respectively, $\zeta$). This implies that
\begin{equation}\label{equ6}
 \int_{\partial \tilde{U}_j\times\partial \tilde{U}_j}
g_1(\chi^{-1}(w))G_{\tilde{h}_j}(w,w')g_2(\chi^{-1}(w'))
  dwdw' =\int_{\partial U_j\times\partial U_j}g_1(\zeta)G(\zeta,\zeta')g_2(\zeta')d\zeta d\zeta',
\end{equation}
where $g_1$ and $g_2$ are mean-zero test-functions whose integrals over $\partial U_j$ vanish. Since $(h,g)_\nabla$ is a  Gaussian random variable with zero mean for each $g\in H(\mathbb{C})$,  the distribution of
$h_j=h|_{U_j}$ is uniquely determined by its covariances. So we conclude from (\ref{equ6})
that $h_j$ and $\tilde{h}_j$
have identical distributions for $j=1,2,\dots,m$.

Next, Let $\tau_j$  and $\tilde{\tau}_j$ be random measures corresponding to
$h_j$ and $\tilde{h}_j$ respectively, as defined in (\ref{equ3}) and (\ref{equ3a}).
Then the equivalence of distributions of $h_j$ and $\tilde{h}_j$ implies
that $\tau_j$ and $\tilde{\tau}_j$
have the same laws.
Hence it follows from (\ref{equ8})
that $p_j(x)\sim\tilde{p}_j(x)$, where $\tilde{p}_j(x)=\tilde{\tau}_j[0,x]/\tilde{\tau}_j[0,1],
x\in [0,1)$, that is, $p_j(x)$ is conformally invariant.
Thus we conclude from (\ref{equ7}) that the distribution of the random homeomorphism $\psi_j$
 is identified with   $\chi\circ \psi_j\circ\chi^{-1}$. This completes the proof of the lemma. $\hfill\square$

\section{The extension of random homeomorphisms} \label{erh}

In this section based on the approach of the Beuring-Ahlfors extension \cite{bea}, we shall describe how to
extend the random homeomorphism $\psi_j:\partial U_j\rightarrow \partial U_j$ constructed in Section \ref{home}
to the multiply connected domain $\Omega\subset\mathbb{S}^2$, and then give a geometric estimate of the corresponding distortion
function in terms of $\psi_j$ and a  estimate for the associated Lehto
integral. Finally, we discuss the uniqueness of the random conformal welding  of $\Omega$ induced by $\psi_j$, which
involves the conformal removability of the boundary of $\Omega$.

\subsection{Construction of the extension mapping} \label{CEM}

For the multiply connected domain $\Omega=\mathbb{C}\setminus \cup_{j=1}^m \bar{U}_j$ as defined before,
without loss of generality
we may assume that \begin{equation}\label{equ9a} \mbox{dist}(\partial U_i,\partial U_j)\geq r e^{4\pi}\end{equation}
for any pair $(i,j), i\neq j$, where $r=\max_{1\leq j\leq m}\{r_j\}$. Otherwise, consider
another multiply connected domain $\tilde{\Omega}=\mathbb{C}\setminus \cup_{j=1}^m \tilde{U}_j$
whose boundary components satisfy (\ref{equ9a}), where $\partial U_k$
are replaced by $\partial \tilde{U}_k$
and $r$ is replaced by $\tilde{r}$.  By Koebe's Theorem and analytic
continuations on $\partial U_j$ we can find a conformal mapping $\varphi$ from a domain $N_{\Omega}\supset\Omega$
onto another domain $N_{\tilde{\Omega}}\supset\tilde{\Omega}$ such that $\varphi(\partial U_j)=\partial \tilde{U}_j$
for $j=1,2,\dots,m$. It follows from Lemma \ref{lem3}(ii) that the random homeomorphism
$\psi_j$ on $\partial U_j$ is equivalent in law to the corresponding one on $\partial \tilde{U}_j$.

Write $$N_{U_j}=\{z\in\bar{\Omega}:|z-a_j|<r_je^{4\pi}\}, \ N_{U_j}^c=\bar{\Omega}
\setminus N_{U_j}=\{z\in\bar{\Omega}:|z-a_j|\geq r_je^{4\pi}\}$$
 for $j=1,2\dots,m$, which will be used below.

We first describe how to extend $\psi_1:\partial U_1\rightarrow \partial U_1$ to $\Omega$.
The definition of $\psi_1$ (see (\ref{equ7}) and (\ref{equ8}))
gives that the random homeomorphism
$p_1:\mathbb{R}\rightarrow\mathbb{R}$ satisfies
\begin{equation}\label{equ9}p_1(x+1)=p_1(x)+1,~ p_1(0)=0.
\end{equation}
Thus we can extend $p_1$ to the upper half plane
$\mathbb{H}$ according to the techniques of  Beuring-Ahlfors \cite{bea}.
To be more concrete, we let
\begin{equation}\label{equ10} F_1(x+iy)=\frac{1}{2}\int_0^1(p_1(x+ty)+p_1(x-ty))dt+i\int_0^1(p_1(x+ty)-p_1(x-ty))dt
\end{equation}
for $0<y<1$. Then it is easy to see that $F_1=p_1$ on $\mathbb{R}$ and $F_1$ is a continuously differentiable
homeomorphism.
Furthermore, from (\ref{equ9}) and (\ref{equ10}) we may set $F_1(z)=z+(2-y)M_0$ for $1\leq y\leq 2$, where
$M_0=\int_0^1 p_1(t)dt-1/2$. This implies that we are able to take $F_1(z)\equiv z$ for $y\geq 2$.
In addition, it is clear that one has
\begin{equation}\label{equ10c}F_1(z+k)=F_1(z)+k
\end{equation}
for any $k\in\mathbb{Z}$. Therefore,
$F_1$ is the desired extension of $p_1$ to $\mathbb{H}$.

Hence, the extension of $\psi_1$ to the disk $U_1$, denoted by $\tilde{\Phi}_1$, can be given by
\begin{equation}\label{equ10c1}
\tilde{\Phi}_1(z)=a_1+r_1\exp(2\pi iF_1(\frac{\log ((z-a_1)/r_1)}{2\pi i})),~ z\in U_1.
\end{equation}
It is easy to see from (\ref{equ7}) and (\ref{equ10c}) that $\tilde{\Phi}_1$ is a well-defined
homeomorphism of $U_1$ with
$\tilde{\Phi}_1|_{\partial U_1}=\psi_1$ and $\tilde{\Phi}_1(z)\equiv z$ for
$|z-a_1|\leq r_1e^{-4\pi}$.
Let $U_1^c$ denote the exterior of the disk $U_1$, Then by the reflection around the
circle $\partial U_1$ we obtain
the extension of $\psi_1$ to $U_1^c$, denoted by $\Phi_1^c$,
which is expressed by
\begin{equation}\label{equ10c2}
\Phi_1^c(z):=
a_1+\frac{r_1^2}{\overline{\tilde{\Phi}_1(z^*)}-\bar{a}_1},~ z^*=a_1+\frac{r^2}{\bar{z}-\bar{a}_1},~ z\in U_1^c.
\end{equation}
It is clear that $\Phi_1^c|_{\partial U_1}=\psi_1$.
Moreover, a simple computation gives that
$\Phi_1^c(z)\equiv z\quad  \mbox{for}~ |z-a_1|\geq r_1e^{4\pi}$.
Thus we let
\begin{equation}\label{equ10c3}
\Phi_1=\Phi_1^c|_{\Omega},
\end{equation}
the restriction of $\Phi_1^c$
to $\Omega$.  Then $\Phi_1$ is the extension of $\psi_1$ to $\Omega$
which satisfies $\Phi_1|_{\partial U_1}=\psi_1$ and  $\Phi_1=I$ on $N_{U_1}^c$,
where $I$ denotes the identity mapping.
This, combined with the condition (\ref{equ9a}),
implies that $\Phi_1|_{\partial U_j}=I, j=2,\dots,m$,
where we used the fact that $\partial U_m=\partial U_m^c$.

Secondly, find the extension of $\psi_j$ to $\Omega$ for $j=2,\dots,m-1$.
After $\psi_1$ has been extended to $\Omega$, the homeomorphism  $\psi_2$ on $\partial U_2$
is transformed to
$\Phi_1\circ\psi_2\circ\Phi_1^{-1}$
on $\Phi_1(\partial U_2)$.  Notice that $\Phi_1$ is an identity mapping on $\partial {U}_2$.
So we get that $\Phi_1\circ\psi_2\circ\Phi_1^{-1}=\psi_2$.
Thus applying the same
method to $\psi_2$,  we obtain the extension mapping $\Phi_2$ of $\Omega$
which satisfies $\Phi_2\circ\Phi_1|_{\partial
U_j}=\psi_j, j=1,2$ and $\Phi_2\circ\Phi_1=I$ on $N_{U_1}^c\cap N_{U_2}^c$. It follows from
(\ref{equ9a}) that $\Phi_2\circ\Phi_1|_{\partial U_j}=I, j=3,\dots,m$.
Repeating the above procedure until $\psi_{m-1}$ has been considered, we obtain
$m-1$ extension mappings $\Phi_j$ of $\Omega$ corresponding to $\psi_j$ which
satisfy  $\Phi_{m-1}\circ\dots\circ\Phi_2\circ\Phi_1|_{\partial
U_j}=\psi_j, j=1,2,\dots,m-1$ and $\Phi_{m-1}\circ\dots\circ\Phi_2\circ\Phi_1=I$ on $
N_{U_1}^c\cap N_{U_2}^c\cap\dots\cap N_{U_{m-1}}^c$. Based on the same reason as above, we have
$\Phi_{m-1}\circ\dots\circ\Phi_2\circ\Phi_1|_{\partial U_m}=I$.

Thirdly, seek the extension of $\psi_m$ to $\Omega$. After the
extension mapping $\Phi_j$ of $\psi_j$ have been obtained as above for $j=1,\dots,m-1$,
the homeomorphism $\psi_m$ on $\partial U_m$ is transformed to
$\Phi_{m-1}\dots\circ\Phi_1\circ\psi_m\circ\Phi_1^{-1}\circ\dots,\Phi_{m-1}^{-1}$ on
$\Phi_{m-1}\dots\circ\Phi_1(\partial U_m)$.
Since  $\Phi_{m-1}\circ\dots\circ\Phi_2\circ\Phi_1|_{\partial U_m}=I$,
$\Phi_{m-1}\dots\circ\Phi_1\circ\psi_m\circ\Phi_1^{-1}\circ\dots,\Phi_{m-1}^{-1}=\psi_m$.
Similarly to construction of $\tilde{\Phi}_1$ instead of $\Phi_1^c$,
by the condition (\ref{equ9a}) we can find the extension $\Phi_m$
of $\psi_m$ to $U_m^c$ such
that $\Phi_{m}\circ\dots\circ\Phi_2\circ\Phi_1|_{\partial
U_j}=\psi_j, j=1,2,\dots,m$ and $\Phi_{m}\circ\dots\circ\Phi_2\circ\Phi_1(\Omega)=\Omega$.

Finally, set
$\Phi=\Phi_m\circ\Phi_{m-1}\circ\dots\circ\Phi_2\circ\Phi_1$. Then $\Phi:\Omega\rightarrow\Omega$
is a well-defined random
homeomorphism  which satisfies $\Phi|_{\partial U_j}=\psi_j$
for $j=1,2,\dots,m$. So $\Phi$ is the desired extension of $(\psi_1,\psi_2,\dots,\psi_m)$
to the multiply connected region $\Omega$.

\subsection{Estimates of the distortion function}

Let $K_\Phi$ and $K_{\Phi_j}$ denote the distortion functions of $\Phi$ and $\Phi_j$,
respectively. It follows from the distortion
properties of quasiconformal mappings that
\begin{equation}\label{equ10d1}K_\Phi(z)=K_{\Phi_m}(z_{m-1})
\circ\dots\circ K_{\Phi_2}(z_1)\circ K_{\Phi_1}(z),~ z\in \Omega,
\end{equation}
where $z_j=\Phi_j(z_{j-1}), j=1,2,\dots,m-1$ and $z_0=z$.
Note that $\Phi_j=I$ on $ N_{U_j}^c$ for $j=1,2,\dots,m$,
where $N_{U_j}^c$ is defined in Section \ref{CEM}. So we deduce
from (\ref{equ10d1}) that
\begin{equation}\label{equ10d} K_\Phi(z)=
\begin{cases} K_{\Phi_j}(z),& \text{if $z\in N_{U_j}$}, j=1,2,\dots,m,\\
1, &\text{if $z\in \cap_{j=1}^m N_{U_j}^c$},
\end{cases}
\end{equation}
which reduces all estimates of $K_\Phi$ to  $K_{\Phi_j}|_{N_{U_j}}$
for $j=1,2,\dots,m$. On the other hand, since the distortion properties
are conformally invariant,  we conclude from the construction of $\Phi_j$
that
\begin{equation}\label{equ10d2}K_{\Phi_j}|_{N_{U_j}}= K_{F_j}|_S
\end{equation}
for $j=1,2,\dots,m$,
where $F_j$ is defined in the same way as $F_1$ in Section \ref{CEM},
and $S=\mathbb{R}\times [0,2]$.

We will give a geometric estimate  for the distortion function $K_\Phi$
in terms of the random homeomorphisms on the boundary components of $\Omega$, and estimate the geometric
distortion of an annulus in $\Omega$ under $\Phi$.
To this end, we need to introduce the following notation which is similar to those in \cite{ajks}.
Let $\mathcal{B}_n$ denote the set of all dyadic intervals of length $2^{-n}$, that is,
\[\mathcal{B}_n=\{[k2^{-n},(k+1)2^{-n}]:k\in\mathbb{Z}\},\]
and set $\mathcal{B}=\{\mathcal{B}_n:n\geq 0\}$. For a pair of intervals $\textbf{J}=\{J_1,J_2\}\subset\mathcal{B}$,
let
\[\delta_{\tau_j}(\textbf{J})=\frac{\tau_j(J_1)}{\tau_j(J_2)}+\frac{\tau_j(J_2)}{\tau_j(J_1)},\]
where $\tau_j$ is the random measure \eqref{equ3}. Set
$C_I=\{(x,y): x\in I, 2^{-n-1}\leq y\leq 2^{-n}\}$ for any $I\in\mathcal{B}_n,n>0$, and
$C_I=I\times [1/2,2]$ for $I\in \mathcal{B}_0$. Then $\{C_I\}_{I\in \mathcal{B}}$
paves the strip $S=\mathbb{R}\times[0,2]$. Moreover, for a dyadic interval $I\in \mathcal{B}_n$ we let
$j(I)$ denote the union of $I$ and its neighbors in $\mathcal{B}_n$. Write
\[\Lambda(I):=\{\textbf{J}=(J_1,J_2):J_i\in \mathcal{B}_{n+5}~ \mbox{and}~J_i\in j(I)\}.\]
We define
\[K_{\tau_j}(I):=\sum_{\textbf{J}\in\Lambda(I)}\delta_{\tau_j}(\textbf{J})\]
for $j=1,2,\dots,m$.

In addition, let $B(z,r,R)=\{w\in\mathbb{C}:r<|w-z|<R\}\subset\Omega$
where $0\leq r\leq R<\infty$. The Lehto integral of $K_\Phi$ corresponding to $B(z,r,R)$ is given by
\begin{equation}\label{equ10g} L(z,r,R):=L_{K_\Phi}(z,r,R):=\int_r^R\frac{1}{\int_0^{2\pi}
K_\Phi(z+\rho e^{i\theta})d\theta}\frac{d\rho}{\rho};
\end{equation}
also see \cite{le,aim}. For any  bounded  topological annulus $\tilde{B}\in\mathbb{C}$,
let $D_o(\tilde{B})$ and $D_I(\tilde{B})$ denote its outer diameter and inner diameter, respectively.

\begin{lemma}
\label{lem4}  Let $\Phi$ be the extension of the random homeomorphisms $\psi_1,\psi_2,\dots,\psi_m$
to $\Omega$
as constructed above, and $K_\Phi$ be its distortion function.  Then (a) there exists a constant $M>0$ such that
\begin{equation}\label{equ10e}
 K_\Phi(z)\leq\max\limits_{1\leq j\leq m}\{\sup_{w\in C_I}K_{F_j}(w)\}\leq M\max\limits_{1\leq j\leq m}\{K_{\tau_j}(I)\}
\end{equation}
for each $z\in N_{U_j}\subset\Omega$, where $C_I\subset S\subset\mathbb{H}$ contains the point corresponding to $z$
 via the relationship between $\Phi_j$ and $F_j$ which is given by (\ref{equ10c1}), (\ref{equ10c2}) and (\ref{equ10c3})
(the subindex $1$ is replaced by $j$) for $j=1,2,\dots,m-1$, and by (\ref{equ10c1}) (the subindex $1$ is replaced by $m$) for $j=m$.

(b) It holds that
\begin{equation}\label{equ10f}\frac{D_o(\Phi(B(z,r,R))}{D_I(\Phi(B(z,r,R))}\geq
\frac{1}{16}e^{2\pi^2L_{K_{\Phi_j}}(z,r,R)}
\end{equation}
for any annulus $B(z,r,R)\subset N_{U_j}\subset\Omega ~(j=1,2,\dots,m)$.
\end{lemma}
\noindent\textbf{Proof.} By the relationship between $\Phi_j$ and $F_j$,
we get that for each $z\in N_{U_j}\subset\Omega$, there exists a $w\in C_I\subset S$
corresponding to $z$.  Note that each $F_j$ is a quasiconformal mapping from
$\mathbb{H}$ onto itself, which is obtained by the extension of
the random homeomorphism $p_j:
\mathbb{R}\rightarrow\mathbb{R}$ to $\mathbb{H}$.
So we deduce the first inequality of (\ref{equ10e}) from (\ref{equ10d}) and (\ref{equ10d2}).
Applying \cite[Theorem 2.6]{ajks} to any $F_j$, we have
\[\sup_{w\in C_I}K_{F_j}(w)\leq M_j K_{\tau_j}(I),\]
where $M_j>0$ is a universal constant.  Thus, setting $M=\max_{1\leq j\leq m}M_j$, we conclude
that the second inequality of (\ref{equ10e}) holds.

It follows from the construction of $\Phi$ that
$\Phi(B(z,r,R))=\Phi_j(B(z,r,R))$ for each
annulus $B(z,r,R)\subset N_{U_j}\subset\Omega$.
For the quasiconformal mapping $\Phi_j$, we get from \cite[Lemma 2.3]{ajks} that
\[\frac{D_o(\Phi_j(B(z,r,R))}{D_I(\Phi_j(B(z,r,R))}\geq \frac{1}{16}e^{2\pi^2L_{K_{\Phi_j}}(z,r,R)}.\]
This yields that (\ref{equ10f}) holds.
Thus we finish the proof of the lemma. $\hfill\square$

\

Now, in the light of (\ref{equ10e}) we may define $K_{\tau_j}$ in the upper half-plane $\mathbb{H}$ by  setting
\[ K_{\tau_j}(z):=K_{\tau_j}(I)\]
for $z\in C_I$. Then a lower bound for the Lehto integral (\ref{equ10g}) can be obtained through replacing
$K_\Phi$ by $K_{\tau_j}$. In the same manner, we may define $K_{\nu_j}(z)(z\in\mathbb{H})$ via the modified
Bearling-Ahlfors extension of the periodic homeomorphism produced by the measure $\nu_j$, as in Section \ref{CEM}.
Thus, in order to show that almost surely there exist infinitely many annuli
around each point on $\partial \Omega$ which are not distorted much by the quasiconformal mapping $\Phi$,
we need the following proposition on probabilistic estimates for Lehto
integrals and the almost sure integrability of the distortions.

\begin{prop}
\label{prop1} Let $\beta_j<\sqrt{2}$, and let $K_{\nu_j}$
be defined as above for $j=1,2,\dots,m$.
Then (i) for each $w\in\mathbb{R}$, there exist $\sigma>0$, $r_0>0$ and $\delta(r)>0$
such that for any  positive $r<r_0$
and $\delta<\delta(r)$ the Lehto integral of $K_{\nu_j}$ satisfies the probabilistic estimate
\begin{equation}\label{equ10a}
\mathbb{P}(L_{K_{\nu_j}}(w,r^n,2r)<n\delta)\leq r^{(1+\sigma)n},
n\in\mathbb{N},
\end{equation}
for $j=1,2,\dots,m$, where $L_{K_{\nu_j}}(w,r^n,2r)$ is defined in (\ref{equ10g}).

(ii) Almost surely $K_{\nu_j}\in L^1([0,1]\times[0,2])$ for $j=1,2,\dots,m$.
\end{prop}

\noindent\textbf{Proof.} Notice that the probability law of  $\nu_j$ is equal to
that of $\nu$ in \cite{ajks}
for $j=1,2,\dots,m$.
So for each $K_{\nu_j}$, we deduce from \cite[Theorem 4.1]{ajks} that
there exists $\sigma_j>0$, $r_{j,0}>0$ and  $\delta_j(r)>0$ such that
for any positive $r<r_{j,0}(r=2^{-p},p\in\mathbb{N})$
and $\delta<\delta_j(r)$ the Lehto integral satisfies the estimate
\[\mathbb{P}(L_{K_{\nu_j}}(w,r^n,2r)<n\delta)\leq r^{(1+\sigma_j)n},
 n\in\mathbb{N}.\]
Take $\sigma=\min_{1\leq j\leq m} \sigma_j$,  $r_0=\min(1, \min_{1\leq j\leq m} r_{j,0})$
and $\delta(r)=\min_{1\leq j\leq m}\delta_j(r)$.
Then we can obtain that (\ref{equ10a}) holds.
It is easy to see that (ii) follows from \cite[Lemma 4.5]{ajks}.
This completes the proof of the
proposition. $\hfill\square$

\subsection{Uniqueness of the  welding}

In order to prove the uniqueness of random conformal welding
for the multiply connected domain $\Omega$, we need the following lemma
involving the conformal removability of boundary of
the multiply connected region, which generalizes the conformal removability result \cite{js} of
boundary of the simply connected domain to the multiply connected case.

Recall a compact set $E\subset D$ is confromally removable inside
a domain $D\subset\mathbb{S}^2$, if any homeomorphism of $D$, which is
conformal on $D\setminus E$, is conformal on $D$. Let $\Omega=\mathbb{S}^2\setminus\cup_{j=1}^m \bar{U}_j$
be defined as before, and $\tilde{\Omega}=\mathbb{S}^2\setminus \cup_{j=1}^m \bar{\Omega}_j$ be
any $m$-connected domain in $\mathbb{S}^2$. Koebe's Theorem gives that there exists a conformal mapping $\Psi$  of $\Omega$
onto $\tilde{\Omega}$. This conformal mapping $\Psi$ is called the Koebe mapping.
The uniqueness of conformal welding for  $\Omega$  is a consequence of the following Lemma.

\begin{lemma}
\label{lem5} Let $\tilde{\Omega}=\mathbb{S}^2\setminus \cup_{j=1}^m \bar{\Omega}_j$
 be an $m$-connected region such that
the Koebe mapping $\Psi:\Omega=\mathbb{S}^2\setminus\cup_{j=1}^m \bar{U}_j\rightarrow\tilde{\Omega}$
is $\alpha$-H\"{o}lder continuous
for some $\alpha>0$, where $U_j\subset\mathbb{S}^2$, $j=1,2,\dots,m$, are mutually  disjoint  disks, and $\Omega_j\subset\mathbb{S}^2$, $j=1,2,\dots,m$, are mutually disjoint simply connected domains.
Then the boundary $\partial \tilde{\Omega}=\cup_{j=1}^m \partial\Omega_j$ is conformally removable.
\end{lemma}

\noindent\textbf{Proof.} We first show that every $\partial \Omega_j$ is conformally
removable. It follows from \cite{doy,mar} that the mapping $\Psi$ can be written as
a composition of $m$ conformal mappings of simply connected domains, i.e., there exist conformal mappings
$\Psi_j$, $j=1,2,\dots,m$ of simply connected domains into $\mathbb{S}^2$ such that
\begin{equation}\label{equ10b}
\Psi=\Psi_m\circ\cdots\circ\Psi_2\circ\Psi_1\quad  \mbox{on}~~\Omega.
\end{equation}
Set $U_j^c:=\mathbb{S}^2\setminus U_j$, $j=1,2,\dots,m$. So ${U_j^c}'s$ are
simply connected domains with pairwise disjoint complements, and
$\partial \Omega=\cup_{j=1}^m \partial U_j^c=\cup_{j=1}^m \partial U_j$.
From (\ref{equ10b}) we deduce that $\Psi_m\circ\cdots\circ\Psi_2\circ\Psi_1$
maps $U_j^c$ conformally onto $\mathbb{S}^2\setminus \Omega_j$
for $j=1,2,\dots,m$. Moreover,
$\Psi_m\circ\cdots\circ\Psi_2\circ\Psi_1|_{\partial U_j^c}=\Psi|_{\partial U_j^c}$.
Thus by the assumption of $\Psi$ one sees that
 $\Psi_m\circ\cdots\circ\Psi_2\circ\Psi_1|_{\partial U_j^c}$ is
$\alpha$-H\"{o}lder continuous for any $j$.
Hence,  from the result of conformal removability for the boundary of simply connected domain (see \cite{js}
or\cite[Theorem 2.4]{ajks}),
we get that $\partial \Omega_j=\partial(\mathbb{S}^2\setminus \Omega_j)$
is conformally removable for each $j=1,\cdots, m$.

Next, we demonstrate that $\partial \tilde{\Omega}$ is conformally
removable. Indeed, let $\chi$ be any homeomorphism of $\mathbb{S}^2$ which
is conformal off $\partial \tilde{\Omega}=\cup_{j=1}^m \partial \Omega_j$.
Since $\partial\Omega_j\subset\mathbb{S}^2$, $j=1,\cdots,m$,
are mutually disjoint compact sets,  there
exists a simply connected neighborhood $N_{\partial\Omega_j}\supset\partial\Omega_j$ for each $\partial\Omega_j$
such that $N_{\partial\Omega_j}$, $j=1,2,\dots,m$, are pairwise disjoint.
Thus $\chi|_{N_{\partial\Omega_j}}$ is
a homeomorphism of $N_{\partial\Omega_j}$ which is conformal off $\partial \Omega_j$.
Since we have proved that each $\partial\Omega_j$ is conformally removable, $\chi|_{N_{\partial\Omega_j}}$
may be extended conformally to the whole neighborhood $N_{\partial\Omega_j}$. This yields that
$\chi$ is conformal in the whole sphere $\mathbb{S}^2$.  So we get from the definition
of conformal removability that
$\partial \tilde{\Omega}=\cup_{j=1}^m \Omega_j$ is conformally removable.  This completes
the proof of the lemma. $\hfill\square$

\begin{prop}
\label{prop2} Let $\Omega\subset\mathbb{S}^2$ be an  $m$-connected domain
whose complement is a union of mutually disjoint closed disks $\bar{U}_j, j=1,2,\dots,m$
for any fixed integer $m\geq 2$, and let $\psi_j:\partial U_j\rightarrow \partial U_j$ be
random homeomorphisms on the boundary components $\partial U_j$ of $\Omega$. Suppose that there exists a random conformal mapping $f$ from $\Omega$ into
$\mathbb{S}^2$ and a random conformal mapping $g_j$ from $U_j$ into $\mathbb{S}^2$ such that
their boundary values satisfy
(\ref{eq2}).
Assume further that $f$ is  $\alpha$-H\"{o}lder continuous on the boundary $\partial \Omega=\cup_{j=1}^m\partial
U_j$ (equivalently $g_j$ is $\alpha$-H\"{o}lder continuous on $\partial U_j$ for each $j$). Then the welding is unique:
any other tuple of welding mappings $(\tilde{f},\tilde{g}_1,\dots,\tilde{g}_m)$ corresponding to $\psi_j$
 is of the form
\begin{equation}\label{equ111}
\tilde{f}=\chi\circ f,\quad \tilde{g}_j=\chi\circ g_j,
\end{equation}
where $\chi:\mathbb{S}^2\rightarrow\mathbb{S}^2$ is a M\"{o}bius transformation.
\end{prop}

\noindent\textbf{Proof.} Assume that $(\tilde{f},\tilde{g}_1,\dots,\tilde{g}_m)$ is another
tuple of welding mappings admitted by $\psi_j$, $j=1,\cdots,m$. Note that the boundary values of
both $(f,g_1,\dots,g_m)$ and $(\tilde{f},\tilde{g}_1,\dots,\tilde{g}_m)$ satisfy the
welding condition (\ref{eq2}).
So we have
\[\tilde{f}^{-1}\circ \tilde{g}_j(z)=\phi_j(z)=f^{-1}\circ g_j(z),\quad z\in \partial U_j\]
for $j=1,2,\dots,m$.
This, combined with the equality $f(\Omega)\cup f(\partial \Omega)\cup_{j=1}^m g_j(U_j)=
\tilde{f}(\Omega)\cup \tilde{f}(\partial \Omega)\cup_{j=1}^m \tilde{g}_j(U_j)=\mathbb{S}^2$,
implies that
\begin{equation*} \chi(z)=
\begin{cases} \tilde{f}\circ f^{-1}(z),& \text{if $z\in f(\Omega)$},\\
\tilde{g}_j\circ g_j^{-1}(z), &\text{if $z\in g_j(U_j)$}(j=1,2,\dots,m)
\end{cases}
\end{equation*}
defines a homeomorphism of $\mathbb{S}^2$ which is conformal off
$f(\partial\Omega)=\cup_{j=1}^m f_j(\partial U_j)$.
Since $f$ is $\alpha$-H\"{o}lder continuous on the boundary $\partial \Omega=\cup_{j=1}^m\partial
U_j$,  it follows from Lemma \ref{lem5}  that $f(\partial\Omega)$ is conformally removable. The definition
of conformal removability gives that $\chi$ can
be extended conformally to the entire sphere $\mathbb{S}^2$.
Hence we get that $\chi$ is a M\"{o}bius transformation
of $\mathbb{S}^2$ which satisfies (\ref{equ111}). So we finish the proof. $\hfill\square$

\section{Random conformal welding theorems } \label{RCW}

In this section we will establish random conformal welding theorems for the multiply
connected domain $\Omega$. The random welding problem
for $\Omega$ is first reduced to solving the associated Beltrami equation (Theorem \ref{thm1}). Then from
the existence of the solution to the Beltrami equation, and the uniqueness of the solution (Proposition \ref{prop2}),  we obtain
the desired solution to the random welding problem (Theorem \ref{thm2}). Finally, we present one result of random conformal welding
with the random homeomorphism on each boundary component of $\Omega$
 arising from two independent Gaussian free fields (Theorem \ref{thm3}).

\begin{theorem}
\label{thm1} Suppose that $\Omega$ is an $m$-connected domain
of the Riemann sphere $\mathbb{S}^2$ whose complement is a union of disjoint closed disks $\bar{U}_j$, $j=1,\cdots,m$ for any
fixed integer $m\geq 2$. Assume that $\beta_j<\sqrt{2}$, $j=1,2,\dots,m$. Let $\psi_j:\partial U_j\rightarrow\partial U_j$
be the  random circle homeomorphism corresponding to $\beta_j$ as defined in Section \ref{home},
and let $\Phi:\Omega\rightarrow\Omega$
be the homeomorphism which extends $\psi_1,\psi_2,\dots,\phi_m$ as constructed in Section \ref{CEM}. Let
\[\lambda=\lambda_\Phi:=\frac{\partial_{\bar{z}}\Phi}{\partial_z\Phi}\]
be the complex dilatation of the extension on $\Omega$, and set $\lambda=0$ on $\bar{U}_j$ for $j=1,2,\dots,m$.

Then almost surely there exists a random homeomorphic $W_{\mbox{loc}}^{1,1}$-solution $F:
\mathbb{C}\rightarrow\mathbb{C}$ to the Beltrami equation
\begin{equation}\label{equ11}
\partial_{\bar{z}}F=\lambda\partial_z F,\quad \mbox{a.e. in}\quad \mathbb{C},
\end{equation}
which satisfies the normalization $F(z)=z+o(z)$ as $z\rightarrow\infty$. In addition, there exists a positive constant $\alpha$
such that the restriction $F|_{\partial U_j}:\partial U_j\rightarrow\mathbb{C}$ is a.s. $\alpha$-H\"{o}lder continuous
for $j=1,2,\dots,m$.
\end{theorem}

\noindent\textbf{Proof.} In the proof we adopt the idea of \cite[Theorem 5.1]{ajks}
or \cite[Theorem 20.9.4]{aim}. First of all, we derive the estimates for $K_\Phi$ and the Lehto integral of the distortion function $K_{\tau_j}$ in the upper half-plane $\mathbb{H}$.
For each integer $n\geq 1$ we
take $M_n=[r^{-(1+\sigma/2)n}]\in\mathbb{N}$
where $\sigma$ is as in Proposition \ref{prop1}(i). For any $j:1\leq j\leq m$, we let
\[z_{n,k}^{(j)}:=a_j+r_je^{2\pi ik/M_n}\]
for $k=1,2,\dots,M_n$, and write $Z_n^{(j)}:=\{z_{n,1}^{(j)},z_{n,2}^{(j)},\dots,
z_{n,M_n}^{(j)}\}$. Then for each $ j=1,\cdots,m$, the distance from $\partial U_j$ to the set $Z_n^{(j)}$ is bounded by
$\pi r_j/M_n\sim r^{(1+\sigma/2)n}$. For a given $n\geq 1$ and
$k\in\{1,2,\dots,M_n\}$, let $E_{n,k}^{(j)}$ denote the event
\[E_{n,k}^{(j)}=\{\omega:L_{K_{\nu_j}}(k/M_n,r^n,2r)<n\delta\}\]
for $j=1,2,\dots,m$ and write $E_n^{(j)}=\cup_{k=1}^{M_n}E_{n,k}^{(j)}$. Then we obtain from Proposition \ref{prop1}(i) that
\[\sum_{n=1}^\infty\mathbb{P}(E_n^{(j)})\leq\sum_{n=1}^\infty\sum_{k=1}^{M_n}\mathbb{P}(E_{n,k}^{(j)})
\leq\sum_{n=1}^\infty M_n r^{(1+\sigma)n}\leq\sum_{n=1}^\infty r^{\sigma n/2}<\infty\]
for $j=1,2,\dots,m$. Borel-Cantelli's lemma implies that for almost every $\omega$ there exists an $n_0(\omega)\in\mathbb{N}$ such that $\omega$
 belongs to the complement of the event $\cup_{n=n_0}^\infty E_n^{(j)}, j=1,2,\dots,m$.

It follows from Lemma \ref{lem2}(a), combined with the definitions of $K_{\tau_j}$ and $K_{\nu_j}$,
that
\[K_{\tau_j}\leq X_j^2 K_{\nu_j}\]
for $j=1,2,\dots,m$, where almost surely $X_j<\infty$. So by Lemma \ref{lem4} (a) and Proposition \ref{prop1}(ii)
we deduce that almost surely
\[\int_{[0,1]\times [0,2]}K_{F_j}(w)dw\leq M \max_{1\leq j\leq m}\int_{[0,1]\times [0,2]}K_{\tau_j}(w)dw
\leq M \max_{1\leq j\leq m}\{X_j^2\int_{[0,1]\times [0,2]}K_{\nu_j}(w)dw\}<\infty.\]

Thus, for a fixed event $\omega_0$ and its corresponding extension $\Phi$ on $\Omega$ with the complex dilatation $\lambda$,
we get from (\ref{equ10d}) that the distortion
\begin{equation} \label{K}
K=K_\Phi=\frac{1+|\lambda|}{1-|\lambda|}
\end{equation}
satisfies
\[K_\Phi(z)=K_{\Phi_j}(z)=K_{F_j}(w)\leq M K_{\tau_j}(w)\leq M X_j(\omega_0)^2K_{\nu_j}(w)\]
for $z\in N_{U_j}\subset\Omega, j=1,2,\dots,m$, where $w\in\mathbb{H}$ is a point corresponding to
$z$ through the mappings $\Phi_j$ and $F_j$;
and $K_\Phi(z)=1$ for $w\in\Omega\setminus\cup_{j=1}^mN_{U_j}$. Also,
$K_{\nu_j}\in L^1\cap L_{\mbox{loc}}^\infty([0,1]\times(0,2])$, and for any $n\geq n_0(\omega_0)$ and
$k\in\{1,2,\dots, M_n\}$ we obtain from the definition of Lehto's integral that
\begin{equation}\label{equ11b}
L_{K_{\tau_j}}(k/M_n,r^n,2r)\geq X_j(\omega_0)^{-2}L_{K_{\nu_j}}(k/M_n,r^n,2r)
\geq n\min_{1\leq j\leq m}\{\delta X_j(\omega_0)^{-2}\}=: n\tilde{\delta}.
\end{equation}

Secondly, we consider the sequence $\{\lambda_l\}$ whose limit is $\lambda$ and
show that the solutions
to Beltrami equations with $\lambda_l$ converge uniformly on compact sets of $\Omega$ to
 the solution of (\ref{equ11}) in the
light of Arzela-Ascoli's thorem.
To this end, we take
\[\lambda_l:=\frac{l}{l+1} \lambda,\quad l\in\mathbb{N}.\]
Let $F_l$ denote the corresponding solution of the Beltrami equation with coefficient
$\lambda_l$, which satisfies the normalization
\begin{equation}\label{equ11b1}
F_l=z+o(1)\quad \mbox{as}\quad z\rightarrow 0.
\end{equation}
Then each
$F_l$ is a quasiconformal homeomorphism of $\mathbb{C}$.

Let $G_l$ denote the inverse mapping of $F_l$, i.e., $G_l=F_l^{-1}$.
Note that $\Omega=U_m^c\setminus \cup_{j=1}^{m-1}U_j$,
where $U_m^c=\{z\in\mathbb{C}:|z-a_m|<r_m\}$
is a disk with finite radius $r_m$, and $\lambda_l(z)\equiv 0$ for $z\in U_j, j=1,2,\dots,m-1$.
According to \cite[Lemma 20.2.3]{aim},
we can deduce that
\begin{equation}\label{equ11a}|G_l(w_1)-G_l(w_2)|\leq \frac{16\pi^2}{\log(e+|w_1-w_2|^{-1})}(|w_1|^2+|w_2|^2+\int_{\Omega}
\frac{1+|\lambda_l(z)|}{1-|\lambda_l(z)|}dz)
\end{equation}
for any $ w_1,w_2\in\mathbb{C}$. At the same time, observe that for $z\in N_{U_j}(j=1,2,\dots,m)$, one has
\[\frac{1+|\lambda_l(z)|}{1-|\lambda_l(z)|}\leq K_\Phi(z)\leq M K_{\tau_j}(w),\]
where $w$ is the point in $\mathbb{H}$ corresponding to $z$ through $\Phi_j$
and $F_j$, and $K_{\tau_j}\in L^1([0,1]\times [0,2])$.
This gives that the integral in (\ref{equ11a})
is uniformly bounded with respect to $l$.  Hence we conclude that  for any $l\in\mathbb{N}$,
the left hand side of the inequality in (\ref{equ11a}) tends to zero as $|w_1-w_2|\rightarrow 0$, which yields
that the sequence of $\{G_l\}$ forms an equicontinuous family.

We next show that the family $\{F_l\}$ is equicontinuous, too. For any $z\in \Omega$
we set $d=\min_{1\leq j\leq m}
\mbox{dist} (z,\partial U_j)/2$.
It is easy to see that $K$ in \eqref{K} is bounded on $B(z,d)$. Since
\[K_l:=K_{F_l}(\cdot)\leq K\]
for any $l\geq 1$, we get that for $b\in(0,d/2)$
\[L_{K_l}(z,b,r_m)\geq L_K(z,b,d)\geq \frac{1}{\|K\|_{L^\infty(B(w,d))}}\log\frac{d}{b}
\rightarrow 0\]
as $b\rightarrow 0$. In addition, we get from Koebe's theorem or \cite[Corollary 2.10.2]{aim} that
\[F_l(2U_m^c)\subset 5U_m^c,\]
which implies that $\mbox{diam}(F_l(B(z,r_m))\leq 5r_m$. Thus by Lemma \ref{lem4}(b) we can deduce that
$\mbox{diam}(F_l(B(z,b)))$ converges to $0$ uniformly in $l$, as $b\rightarrow 0$.
This gives that $F_l$ is  equicontinuous at every point $z\in\Omega$.
Since $F_l$ is conformal on $U_j$ for $j=1,2,\dots,m$ and satisfies (\ref{equ11b1}),
the equicontinuity of $F_l$ in $U_j$ follows from Koebe's theorem.

Now we will prove the equicontinuity of $F_l$ on $\partial \Omega=\cup_{j=1}^m\partial U_j$. It suffices to
prove local equicontinuity  on points of $[0,1]$ for the families
\[\tilde{F}_l^{(j)}(t)=F_l(a_j+r_je^{2\pi it}),~j=1,2\dots,m.\]
Assume that $n\geq n_0(\omega_0)$. Observe that $\mbox{diam}(\tilde{F}_l^{(j)}(B(k/M_n,2r)))\leq
\mbox{diam}(F_l(B(z_{n,k}^{(j)},r_m)))\leq 5r_m$, which, combined with Lemma \ref{lem4}(b) and (\ref{equ11b}), implies that
\begin{equation}\label{equ12}
\mbox{diam}(\tilde{F}_l^{(j)}(B(k/M_n,r^n)))\leq\mbox{diam}
(\tilde{F}_l^{(j)}(B(k/M_n,2r)))16e^{-2\pi^2n\tilde{\delta}}
\leq 80r_me^{-n\tilde{c}}.
\end{equation}
Since the set $Z_n^{(j)}=\{z_{n,1}^{(j)},z_{n,2}^{(j)},\dots,
z_{n,M_n}^{(j)}\}$ is evenly spread on $\partial U_j$ for each $j\in\{1,2,\dots,m\}$,
the balls $B(z_{n,j}^{(j)},r^{n+1})$ cover the $r^{n+2}$-neighborhood of $\partial U_j$ in such a way that any two points in this neighborhood, whose distance is less than or equal to $r^{n+2}$,
lie in the same ball.  Note that this holds for any $n\geq n_0(\omega_0)$. So we can deduce from (\ref{equ12})
that there are $\epsilon_0>0$ and $\alpha>0$ such that, uniformly in $l$,
\begin{equation}\label{equ13}
|F_l(\tilde{z})-F_l(z)|\leq C|\tilde{z}-z|^\alpha
\end{equation}
when $|\tilde{z}-a_j|=r_j,r_1-\epsilon_0\leq |z-a_j|\leq r_j+\epsilon_0$ and $|z-\tilde{z}|\leq \epsilon_0$.
In fact, we may take $\alpha=\tilde{c}/\log(1/r)$.
This implies that the family $\{F_l\}$ is equicontinous on $\partial \Omega=\cup_{j=1}^m \partial U_j$. Hence we obtain
that $\{F_l\}$ is equicontinous on $\mathbb{S}^2$. Thus applying Arzela-Ascoli's theorem and
passing to a limit,  we obtain a $W^{1,1}$-homeomorphic
solution $F(z)=\lim_{l\rightarrow\infty}F_l(z)$ to the Beltrami equation (\ref{equ11}).

Finally, we get from $(\ref{equ13})$ that $F:\partial U_j\rightarrow\mathbb{C}$ is $\alpha$-H\"{o}lder continuous
for $j=1,2,\dots,m$.
Note that $F$ is analytic in $U_j, j=1,2,\dots,m$, and satisfies (\ref{equ11b1}). So we conclude that $F$ is $\alpha$-H\"{o}lder continuous
on the components $\bar{U}_j$ of $\mathbb{S}^2\setminus \Omega$. This completes the proof. $\hfill\square$

\begin{theorem} \label{thm2}
 Let $\Omega$ be an  $m$-connected domain
of the Riemann sphere $\mathbb{S}^2$ whose complement is a union of disjoint closed disks $\bar{U}_j, j=1,2,\dots,m$
for each fixed integer $m\geq 2$, and suppose that $\psi_j=\psi_{\omega,j}:\partial U_j\rightarrow\partial U_j$ is the random
homeomorphism on the boundary component $\partial U_j$ of $\Omega$ for $0<\beta_j<\sqrt{2}$,
with the exponential GFF as its derivatives, which is defined in Section \ref{home}.
Then almost surely in $\omega$, there exist random conformal mappings
$f:\Omega\rightarrow\mathbb{S}^2$ and
$g_j:U_j\rightarrow\mathbb{S}^2$, respectively, such that their boundary values satisfy
(\ref{eq2}), which produce $m$ mutually disjoint  random Jordan curves $\Gamma_j=\Gamma_{\omega,j}
=f(\partial U_j)=g_j(\partial U_j)$, $j=1,\cdots,m$,
depending on $(\beta_1,\beta_2,\dots,\beta_m)$. Moreover, almost surely in $\omega$, these
Jordan curves $\Gamma_j$ are unique, up to composing with a  M\"{o}bius transformation
$\chi=\chi_\omega$ of the Riemann sphere $\mathbb{S}^2$.
\end{theorem}

\noindent\textbf{Proof.} We first extend $(\psi_1,\psi_2,\dots,\psi_m)$ to a homeomorphism
$\Phi:\Omega\rightarrow\Omega$ as proceeded in Section \ref{CEM} and define a complex dilatation $\lambda(z)$
corresponding to $\Phi$ by
\begin{equation*} \lambda(z):=
\begin{cases} {\partial_{\bar{z}}\Phi}/{\partial_z\Phi},& \text{if $z\in \Omega$},\\
0, &\text{if $z\in \bar{U}_j, j=1,2,\dots,m$}.
\end{cases}
\end{equation*}
Then Theorem \ref{thm1} gives that there must be a homeomorphic solution $F$ to
the Beltrami equation (\ref{equ11}). It is clear that $F$ is conformal in $U_j$.
Thus we put $g_j=F|_{U_j}$ for $j=1,2,\dots,m$. Next, notice that $K_F(z)$ is locally bounded in $\Omega$.
So by the uniqueness of the solution to the Beltrami equation we deduce that there exists
a conformal homeomorphism $f$ defined on $\Omega$ such that
\begin{equation}\label{equ14}F(z)=f\circ \Phi(z),~z\in \Omega.
\end{equation}
Since $\partial U_j$, $j=1,\cdots,m$, are pairwise disjoint and $F$ is homeomorphic on $\mathbb{C}$,
we get that the image boundary
\begin{equation}\label{equ15}\partial f(\Omega)=\cup_{j=1}^m\partial f(U_j)
=\cup_{j=1}^m\partial g(U_j)=F(\cup_{j=1}^m
\partial U_j)\end{equation}
is a union of mutually disjoint  Jordan curves $\Gamma_j$,
where $\Gamma_j=f(\partial U_j)=g_j(\partial U_j)$.
This implies that the conformal mappings $f$ and $g_j$ can be extended to $\partial\Omega=\cup_{j=1}^m
\partial U_j$. Thus from (\ref{equ14}) and the definitions of $g_j$ and $\Phi$, we
deduce easily that $f$ and $g_j$ satisfy  (\ref{eq2})
on the boundary $\cup_{j=1}^m
\partial U_j$ of $\Omega$. Finally, it follows from Theorem \ref{thm1} again that $g_j$ is H\"{o}lder continuous
on $U_j$. This, combined with Proposition \ref{prop2},
implies that the random welding curves $\Gamma_j$ are unique up to
composing with a M\"{o}bius transformation of $\mathbb{S}^2$.
So we finish the proof of the theorem.  $\hfill\square$

\begin{theorem}
\label{thm3} Let $\Omega$ be an  $m$-connected domain
of the Riemann sphere $\mathbb{S}^2$ whose complement is a union of disjoint closed disks $\bar{U}_j, j=1,2,\dots,m$
for each fixed integer $m\geq 2$, and let $0\leq\beta_j^+,\beta_j^-\leq\sqrt{2}$. Suppose that $\psi_j^+$
and $\psi_j^-$ are two independent copies of the random homeomorphism of $\partial U_j$ as defined
in Section \ref{crm}, associated with parameters $\beta_j^+$
and $\beta_j^-$, and two independent GFFs, respectively. Then, almost surely in $\omega$, there exist random conformal mappings
$f:\Omega\rightarrow\mathbb{S}^2$ and
$g_j:U_j\rightarrow\mathbb{S}^2$, respectively, such that their boundary values satisfy
\begin{equation}\label{equ15a}
f^{-1}\circ g_j=\psi_j^+\circ(\psi_j^-)^{-1},
\end{equation}
which produce $m$ mutually disjoint random Jordan curves $\Gamma_j=\Gamma_{\omega,j}=f(\partial U_j)=g_j(\partial U_j)$
depending on $(\beta_1^+,\dots,\beta_m^+;\beta_1^-,\dots,\beta_m^-)$ for $j=1,2,\dots,m$. Moreover,
almost surely in $\omega$, these
Jordan curves $\Gamma_j$ are unique, up to composing with a  M\"{o}bius transformation
$\chi=\chi_\omega$ of the Riemann sphere $\mathbb{S}^2$.
\end{theorem}
\noindent\textbf{Proof.} First,  applying the same method as in Section \ref{CEM}, we extend the boundary homeomorphisms
$\psi_1^+,\psi_2^+,\dots,\psi_m^+$
 to $\Omega$ and denote by $\Phi^+$ the corresponding extension. Moreover, let $\lambda^+$ denote the dilatation of $\Phi^+$
in $\Omega$. Similarly, we may construct the Beurling-Ahlfors extensions of $\psi_j^-$ to $U_j$, $j=1,\cdots,m$,  as in Section \ref{CEM},
and let $\Phi_j^-$ and $\lambda_j^-$ stand for the associated extensions and dilatations respectively. Write
\begin{equation*} \lambda(z)=
\begin{cases} \lambda^+(z),& \text{if $z\in \Omega$},\\
\lambda_j^-, &\text{if $z\in \bar{U}_j, j=1,2,\dots,m$}.
\end{cases}
\end{equation*}
From the specific construction of these extensions and the condition (\ref{equ9a}),
it is easy to see that $\lambda$ has a compact support in $\mathbb{C}$.
Since the estimates for the Lehto integral of the distortion function
\[K(z)=\frac{1+|\lambda|}{1-|\lambda|}\]
in the current situation are equal to those presented in Proposition 1,  carrying through
the same  proof as the one of Theorem \ref{thm1}
with only notational changes we can find as before a solution to the Beltrami  equation
\[\frac{\partial F}{\partial\bar{z}}(z)=\lambda(z)\frac{\partial F}{\partial z}(z)\]
for almost every $z\in\mathbb{C}$, which satisfies the normalization (\ref{equ11b1}).
At the same time, $F|_{\partial U_j}$ is  H\"{o}lder continuous for $j=1,2,\dots,m$.

Next, due to the uniqueness of the solution to the Beltrami equation, there exist conformal
mappings $f:\Omega\rightarrow\mathbb{S}^2$ and $g_j:U_j\rightarrow\mathbb{S}^2$
such that
\begin{equation}\label{equ16}F(z)=f\circ \Phi^+(z),~z\in \Omega
\end{equation}
and
\begin{equation}\label{equ17}F(z)=g_j\circ\Phi_j^-(z),~ z\in U_j
\end{equation}
for $j=1,2,\dots,m$. Hence, arguing as in Theorem 1 we obtain $m$ mutually disjoint  random
Jordan curves $\Gamma_j=\Gamma_{\omega,j}=f(\partial U_j)=g_j(\partial U_j)$ which
depend on $(\beta_1^+,\dots,\beta_m^+;\beta_1^-,\dots,\beta_m^-)$. In the same manner,
from (\ref{equ16}) and (\ref{equ17}) we deduce that $f$ and $g_j$ satisfy
(\ref{equ15a})
on $\partial U_j, j=1,2,\dots,m$, which shows that the mappings $f$ and $g_j$ solve the stated  welding problem.

Finally, note that $F|_{\partial U_j}$ and $(\psi_j^-)^{-1}$ are H\"{o}lder continuous, which implies
$g_j|_{\partial U_j}$, $j=1,\cdots,m$, are  H\"{o}lder continuous, too. So we deduce from Proposition \ref{prop2} that $\Gamma_j$ , $j=1,\cdots,m$, are unique up a  M\"{o}bius transformation
$\chi=\chi_\omega$ of the Riemann sphere $\mathbb{S}^2$. This completes the proof of the theorem. $\hfill\square$

\end{document}